\input amstex
\magnification=1200

\def\o{\omega}

\def\w{\wedge}

\def\va{\varphi}

\def\cen{\centerline}

\cen{\bf H\"older continuity of solutions to the Monge-Amp\`{e}re equations} 
\cen{\bf on compact K\"ahler manifolds} 
\vskip 0,5 cm
\noindent
\cen{PHAM HOANG HIEP}
\vskip 0,5 cm
\noindent
{ABSTRACT.} We study H\"older continuity of solutions to the Monge-Amp\`{e}re equations on compact K\"ahler manifolds. In [DNS] the authors have shown that the measure $\omega_u^n$ is moderate if $u$ is H\"older continuous. We prove a theorem which is a partial converse to this result.
\medskip
\noindent
$\o$-plurisubharmonic functions, compact K\"ahler manifolds.
\vskip 0,5 cm
\noindent
{\bf 1. Introduction}
\medskip
\noindent
Let $X$ be a compact $n$-dimensional K\"aler manifold equipped with a fundamental form $\omega$ satisfying $\int\limits_X\o^n=1$. An upper semicontinuous function $\va :\ X\to [-\infty,+\infty)$ is called $\o$-plurisubharmonic ($\o$-psh) if $\va\in L^1(X)$ and $\o_\varphi :=\o + dd^c\va\geq 0$. By PSH$(X,\o )$ (resp. PSH$^-(X,\o )$) we denote the set of $\o$-psh (resp. negative $\o$-psh) functions on $X$. The complex Monge-Amp\`{e}re equation $\o_u^n = f \o^n$ was solved for smooth positive $f$ in the fundamental work of S. T. Yau (see [Yau]). Later S. Kolodziej showed that there exists a continuous solution if $f\in L^p (\omega^n)$, $f\geq 0$, $p>1$ (see [Ko2]). Recently in [Ko5] he proved that this solution is H\"older continuous in this case (see also [EGZ] for the case $X=\bold CP^n$). In Corollary 1.2 in [DNS] the authors have shown that the measure $\omega_u^n$ is moderate if $u$ is H\"older continuous. The main result is the following theorem which give a partial answer to the converse problem:
\medskip
\noindent
{\bf Theorem A.} {\sl Let $\mu$ be a non-negative Radon measure on $X$ such that 
$$\mu (B(z,r))\leq A r^{ 2n-2+\alpha },$$ 
for all $B(z,r)\subset X$ ($A,\alpha >0$ are constants). Then for every $f\in L^p(d\mu )$ with $p>1$, $\int\limits_X f d\mu =1$, there exists a H\"older continuous $\o$-psh function $u$ such that $\o_u^n = fd\mu $.}
\medskip
\noindent
The following results are simple applications of Theorem A:
\medskip
\noindent
{\bf Corollary B.} {\sl Let $\varphi\in$PSH$(X,\o )$ be a H\"older continuous function. Then for every $f\in L^p(\o_\varphi\w\o^{n-1})$ with $p>1$, $\int\limits_X f \o_\varphi\w\o^{n-1} = 1$, there exists a H\"older continuous $\o$-psh function $u$ such that $\o_u^n = f\o_\varphi\w\o^{n-1}$.}
\medskip
\noindent
{\bf Corollary C.} {\sl Let $S$ be a $C^1$ smooth real hypersurface in $X$ and $V_S$ be the volume measure on $S$. Then for every $f\in L^p(dV_S)$ with $p>1$, $\int\limits_X f dV_S =1$, there exists a H\"older continuous $\o$-psh function $u$ such that $\o_u^n = fdV_S$.}
\medskip
\noindent
{\bf Acknowledgments.} The author is grateful to Slawomir Dinew and Nguyen Quang Dieu for valuable comments. The author is also indebted to the referee for his useful comments that helped to improve the paper.
\vskip 0.5 cm
\noindent
{\bf 2. Preliminaries}
\medskip
\noindent
First we recall some elements of pluripotential theory that will be used throughout the paper. Details can be found in [BT1-2], [Ce1-2], [CK], [CGZ], [De1-3], [Di1-3], [GZ1-2], [H], [Ko1-5], [KoTi], [Si1-2], [Ze1-2].
\medskip
\noindent
{\bf 2.1.} In [Ko2] Ko{\l}odziej introduced the capacity $C_X$ on $X$ by
$$C_X(E):=\sup\{\int\limits_E\o_\va^n:\ \va\in\text{PSH}(X,\o),\ -1\leq\va\leq 0\}$$
for all Borel sets $E\subset X$. 
\medskip
\noindent
{\bf 2.2.} In [GZ1] Guedj and Zeriahi introduced the Alexander capacity $T_X$ on $X$ by
$$T_X(E)=e^{-\sup\limits_X V_{E,X}^*}$$
for all Borel sets $E\subset X$. Here $V_{E,X}^*$ is the global extremal $\omega$-psh function for E defined as the smallest upper semicontinuous majorant of $V_{E,X}$ i.e,
$$V_{E,X}(z) =\sup\{\varphi (z):\ \va\in\text{PSH}(X,\o),\ \va\leq 0\ \text{on}\ E\}.$$
\medskip
\noindent
{\bf 2.3.} The following definition was introduced in [EGZ]: A probability measure $\mu$ on $X$ is said to satisfy the condition $\Cal H(\alpha ,A)$ ($\alpha, A > 0$) if 
$$\mu (K)\leq A C_X(K)^{1+\alpha },$$
for any Borel subset $K$ of $X$. 
\medskip
\noindent
A probability measure $\mu$ on $X$ is said to satisfy the condition $\Cal H(\infty)$ if for any $\alpha >0$ there exist $A(\alpha )>0$ dependent on $\alpha$ such that   
$$\mu (K)\leq A (\alpha ) C_X(K)^{1+\alpha },$$
for any Borel subset $K$ of $X$.
\medskip
\noindent
{\bf 2.4.} The following definition was introduced in [DS]: A measure $\mu$ is said to be moderate if for any open set $U\subset X$, any compact set $K\subset\subset U$ and any compact family $\Cal F$ of plurisubharmonic functions on $U$, there are constants $\alpha >0$ such that 
$$\sup\{\int\limits_K e^{-\alpha\varphi} d\mu :\ \varphi\in\Cal F\}<+\infty.$$ 
\medskip
\noindent
{\bf 2.5.} The following class of $\o$-psh functions was investigated by Guedj and Zeriahi in [GZ2]:
$$\Cal E(X,\o )=\{\va\in\text{PSH}(X,\o ):\ \lim\limits_{j\to\infty }\int\limits_{\{\va >-j\}}\o_{\max (\va ,-j)}^n=\int\limits_X\o^n=1\}.$$
Let us also define
$$\Cal E^- (X,\o )=\Cal E(X,\o )\cap\text {PSH}^-(X,\o ).$$
We refer to [GZ2] for the properties of the class $\Cal E(X,\o)$.
\medskip
\noindent
{\bf 2.6.} $S$ is called a $C^1$ smooth real hypersurface in $X$ if for all $z\in X$ there exists a neighborhood $U$ of $z$ and $\chi\in C^1(U)$ such that $S\cap U=\{z\in U:\ \chi (z)=0\}$ and $D \chi (z)\not = 0$ for all $z\in S\cap U$.
\medskip
\noindent
Next we state a well-known result needed for our work.
\medskip
\noindent
{\bf 2.7. Proposition.} {\sl Let $\mu$ be a non-negative Radon measure on $X$ such that $\mu (B(z,r))\leq A r^{ 2n-2+\alpha }$ for all $B(z,r)\subset X$ ($A,\ \alpha >0$ are constants). Then $\mu\in\Cal H(\infty )$.}
\medskip
\noindent
{\sl Proof.} By Theorem 7.2 in [Ze2] and Proposition 7.1 in [GZ1] we can find $\epsilon, C>0$ such that
$$\mu (K)\leq A h^{2n-2+\alpha }(K)\leq \frac {A C} {\alpha} T_X (K)^{\epsilon\alpha}\leq \frac {A Ce} {\alpha} e^{-\frac {\epsilon\alpha} { C_X (K)^{\frac 1 n } } },$$
for all Borel subsets $K$ of $X$, where $h^{2n-2+\alpha }$ is the Hausdorff content of dimension $2n-2+\alpha$.  This implies that $\mu\in\Cal H(\infty )$.
\vskip 0.5 cm
\noindent
{\bf 3. Stability of the solutions}
\medskip
\noindent
The stability estimate of solutions to the Monge-Amp\`ere equations on compact K\"ahler manifolds was obtained by Kolodziej ([Ko2]). Recently, in [DZ] S. Dinew and Z. Zhang proved a stronger version of this estimate. We will show a generalization of the stability theorem by S. Kolodziej. As a first step we have the following proposition. This proof follows ideas of the proof of Theorem 2.5 in [DH]. We include a proof for the reader's convenience.
\medskip
\noindent
{\bf 3.1. Proposition.} {\sl Let $\varphi,\psi\in\Cal E^-(X,\o)$ be such that $\omega_{\varphi}^n\in\Cal H(\alpha ,A)$. Then there exist constants $t\in\bold R$ and $C(\alpha ,A)\geq 0$ such that} 
$$\int\limits_{\{|\varphi -\psi -t|>a\}}(\o_{\varphi}^n+\o_{\psi}^n)\leq C(\alpha ,A) a^{n+1},$$
here $a=[\int\limits_X ||\o_{\varphi}^n-\o_{\psi}^n||]^{ \frac 1 {2n+3+\frac {n+1}{1+\alpha}}}$.
\medskip
\noindent
{\sl Proof.} Since $\int\limits_{\{|\varphi -\psi -t|>a\}}(\o_{\varphi}^n+\o_{\psi}^n)\leq 2$, it suffices to consider the case when $a$ is small. Set
$$\epsilon =\frac 1 2 \inf \{\int\limits_{\{|\varphi -\psi -t|>a\}}\o_{\varphi}^n:\ t\in\bold R\}$$
Hence
$$\int\limits_{\{|\varphi -\psi -t|\leq a\}}\o_{\varphi}^n\leq 1-2\epsilon$$
for all $t\in\bold R$. Set
$$t_0=\sup\{t\in{\bold R}:\ \int\limits_{\{\varphi<\psi+t+a\}}\o_{\varphi}^n\leq 1-\epsilon\}$$
Replacing $\psi$ by $\psi+t_0$ we can assume that $t_0=0$. Then $\int\limits_{\{\varphi<\psi+a\}}\o_{\varphi}^n\leq 1-\epsilon$ and $\int\limits_{\{\varphi\leq \psi+a\}}\o_{\varphi}^n\geq 1-\epsilon$. Hence 
$$\aligned&\int\limits_{\{\psi<\varphi+a\}}\o_{\varphi}^n=1-\int\limits_{\{\varphi+a\leq \psi\}}\o_{\varphi}^n=1-\int\limits_{\{\varphi\leq \psi+a\}}\o_{\varphi}^n\\
&+\int\limits_{\{\psi- a<\varphi\leq \psi+a\}}\o_{\varphi}^n\leq 1-\epsilon.
\endaligned$$
Since $\int\limits_{\{|\varphi-\psi|\leq a\}}\o_{\varphi}^n\leq 1$ we can choose $s\in [-a+a^{n+2},a-a^{n+2}]$ satisfying $$\int\limits_{\{|\varphi-\psi-s|<a^{n+2}\}}\o_{\varphi}^n\leq 2 {a^{n+1}}.$$
Replacing $\psi$ by $\psi+s$ we can assume that $s=0$. One easily obtains the following inequalities
$$\int\limits_{\{\varphi<\psi+a^{n+2}\}}\o_{\varphi}^n\leq 1-\epsilon,\ \int\limits_{\{\psi<\varphi+a^{n+2}\}}\o_{\varphi}^n\leq 1-\epsilon,\ \int\limits_{\{|\varphi-\psi |<a^{n+2}\}}\o_{\varphi}^n\leq 2 a^{n+1}.\tag 1$$
By [GZ2] we can find $\rho\in \Cal E(X,\omega )$, such that 
$$\omega_{\rho}^n=\frac 1 {1-\epsilon}1_{\{\varphi<\psi\}}\omega_{\varphi}^n+c1_{\{\varphi\geq \psi\}}\omega_{\varphi}^n\ \text{and}\ \sup\limits_X\rho = 0,\tag 2$$ 
($c\geq 0$ is chosen so that the measure has total mass $1$). For simplicity of notation we set $\beta=\frac {n+1} {1+\alpha}$. Set
$$U=\{(1-a^{n+2+\beta }) \varphi < (1-a^{n+2+\beta }) \psi + a^{n+2+\beta }\rho\}\subset\{\varphi<\psi\}.$$
From Theorem 2.1 in [Di3] and (2) we get
$$\omega_{\varphi}^{n-1}\w\omega_{(1-a^{n+2+\beta }) \psi + a^{n+2+\beta }\rho}\geq (1-a^{n+2+\beta })\omega_{\varphi}^{n-1}\w\omega_{\psi}+\frac {a^{n+2+\beta }} {(1-\epsilon)^{\frac 1 n}}\omega_{\varphi}^n,\tag 3$$
on $U$. From Theorem 2.3 in [Di3], Lemma 2.6 in [DH] and (3) we obtain
$$\aligned &(1-a^{n+2+\beta })\int\limits_{U}\omega_{\varphi}^{n-1}\w\omega_{\psi}+\frac {a^{n+2+\beta }} {(1-\epsilon)^{\frac 1 n}}\int\limits_{U}\omega_{\varphi}^n\\
&\leq\int\limits_{U}\omega_{(1-a^{n+2+\beta }) \psi + a^{n+2+\beta }\rho}\w\omega_{\varphi}^{n-1}\\ 
&\leq\int\limits_{U}\omega_{(1-a^{n+2+\beta }) \varphi}\w\omega_{\varphi}^{n-1}=(1-a^{n+2+\beta })\int\limits_{U}\omega_{\varphi}^n+a^{n+2+\beta }\int\limits_{U}\omega\w\omega_{\varphi}^{n-1}\\
&\leq (1-a^{n+2+\beta })(\int\limits_{U}\omega_{\varphi}^{n-1}\w\omega_{\psi}+2a^{2n+3+\beta })+a^{n+2+\beta }\int\limits_{U}\omega\w\omega_{\varphi}^{n-1}.
\endaligned$$
Hence 
$$\frac 1 {(1-\epsilon)^{\frac 1 n}}\int\limits_{U}\omega_{\varphi}^n\leq 2a^{n+1}+\int\limits_{U}\omega\w\omega_{\varphi}^{n-1}.\tag 4$$
From Proposition 3.6 in [GZ1] and (4) we get 
$$\aligned (5)\ \ \ \ &\frac {1} {(1-\epsilon)^{\frac 1 n}}[\int\limits_{\{\varphi\leq \psi-a^{n+2}\}}\omega_{\varphi}^n-C_1(\alpha ,A) a^{n+1}]\\
&\leq\frac {1} {(1-\epsilon)^{\frac 1 n}}[\int\limits_{\{\varphi\leq \psi-a^{n+2}\}}\omega_{\varphi}^n-A [ C_X ( \{\rho\leq -\frac 1 {2 a^{ \beta } }\}) ]^{1+\alpha} ]\\ 
&\leq\frac {1} {(1-\epsilon)^{\frac 1 n}}[\int\limits_{\{\varphi\leq \psi-a^{n+2}\}}\omega_{\varphi}^n-\int\limits_{\{\rho\leq -\frac 1 {2 a^{ \beta } } \}}\omega_{\varphi}^n]\\
&\leq \frac {1} {(1-\epsilon)^{\frac 1 n}}\int\limits_{U}\omega_{\varphi}^n\\
&\leq 2a^{n+1} + \int\limits_{U}\omega\w\omega_{\varphi}^{n-1}\\
&\leq 2a^{n+1} + \int\limits_{\{\varphi<\psi\}}\omega\w\omega_{\varphi}^{n-1},
\endaligned$$
Similarly to $\rho$ we define $\vartheta\in\Cal E(X,\o)$, such that 
$$\omega_{\vartheta}^n=\frac 1 {1-\epsilon}1_{\{\varphi<\psi\}}\omega_{\varphi}^n+l1_{\{\psi\geq \varphi\}}\omega_{\varphi}^n\ \text{and}\ \sup\limits_X\vartheta = 0,$$ 
($l$ plays the same role as $c$ above). Set
$$V=\{ (1-a^{n+2+\beta }) \psi < (1-a^{n+2+\beta }) \varphi + a^{n+2+\beta }\vartheta \}\subset\{\psi<\varphi\}.$$
We get
$$\frac {1} {(1-\epsilon)^{\frac 1 n}}[\int\limits_{\{\psi\leq\varphi -a^{n+2}\}}\omega_{\varphi}^n-C_1(\alpha ,A) a^{n+1}]\leq 2a^{n+1} + \int\limits_{\{\psi<\varphi\}}\omega\w\omega_{\varphi}^{n-1}.\tag 6$$
From (1), (5) and (6) we obtain 
$$\aligned\frac {1} {(1-\epsilon)^{\frac 1 n}}[1-2a^{n+1}-2C_1(\alpha ,A) a^{n+1}]&\leq \frac {1} {(1-\epsilon)^{\frac 1 n}}[\int\limits_{\{|\varphi-\psi|\geq a^{n+1}\}}\omega_{\varphi}^n-2C_1(\alpha ,A) a^{1+\alpha}]\\
&\leq 4 a^{n+1} + 1.
\endaligned$$
Hence
$$\epsilon\leq 1-[ \frac {1-2(C_1(\alpha ,A)+1) a^{n+1}} {4 a^{n+1} + 1} ]^n\leq C_2(\alpha ,A)a^{n+1}.$$ 
This implies that there exists $t\in\bold R$ satisfying
$$\int\limits_{\{|\varphi -\psi -t|>a\}}\o_{\varphi}^n\leq 2C_2(\alpha ,A)a^{n+1}.$$
Finally we have
$$\aligned\int\limits_{\{|\varphi -\psi -t|>a\}}(\o_{\varphi}^n+\o_{\psi}^n)&=2\int\limits_{\{|\varphi -\psi -t|>a\}}\o_{\varphi}^n+\int\limits_{\{|\varphi -\psi -t|>a\}}(\o_{\psi}^n-\o_{\varphi}^n)\\
&\leq 2C_2(\alpha ,A) a^{n+1}+a^{2n+3+\beta }\leq C(\alpha ,A) a^{n+1}.
\endaligned$$
The second step in proving our stability therem is the the following

\noindent
{\bf 3.2. Proposition.} {\sl Let $\varphi,\psi\in\Cal E^-(X,\o)$ be such that $\omega_{\varphi}^n,\omega_\psi^n\in\Cal H(\alpha ,A)$. Then there exist constants $t\in\bold R$ and $C(\alpha ,A)\geq 0$ such that} 
$$C_X(\{ |\varphi -\psi -t|>a \})\leq C(\alpha ,A)a,$$
here $a=[\int\limits_X ||\o_{\varphi}^n-\o_{\psi}^n||]^{ \frac 1 {2n+3+\frac {n+1}{1+\alpha}}}$.
\medskip
\noindent
{\sl Proof.} Since $C_X(\{ |\varphi -\psi -t|>a \})\leq C_X (X)=1$, it suffices to consider the case when $a$ is small. Without loss of generality we can assume that $\sup\limits_X\varphi =\sup\limits_X\psi =0$. By Remark 2.5 in [EGZ] there exists $M(\alpha ,A)>0$ such that $||\varphi||_{L^\infty (X)} < M (\alpha ,A)$, $||\psi||_{L^\infty (X)} < M (\alpha ,A)$. By Proposition 3.1 we can find $t>0$ such that
$$\int\limits_{\{|\varphi -\psi -t|>a\}}(\o_{\varphi}^n+\o_{\psi}^n)\leq C_1(\alpha ,A) a^{n+1}.$$
We consider the case $a<\min (1, \frac 1 {C_1(\alpha ,A)})$. Since $\int\limits_{\{|\varphi -\psi -t|>a\}}(\o_{\varphi}^n+\o_{\psi}^n)<1$ we get $\{|\varphi -\psi -t|>a\}\not = X$. This implies that $|t|\leq\sup\limits_X|\varphi -\psi |+1\leq M(\alpha ,A)+1$. Replacing $\psi$ by $\psi+t$ we can assume that $t=0$ and $||\psi||_{L^\infty (X)} < 2M (\alpha ,A)+1$. Using Lemma 2.3 in [EGZ] for $s=\frac a 2$, $t=\frac a {2(2M (\alpha ,A)+1)}$ we get
$$\aligned C_X (\{ \varphi -\psi <-a\})&\leq C_X(\{ \varphi-\psi<-\frac a 2 - \frac a {2(2M (\alpha ,A)+1)} \})\\
&\leq\frac {2^n(2M(\alpha ,A)+1)^n}{a^n}\int\limits_{\{\varphi -\psi<-a\}}\o_{\varphi}^n\\
&\leq 2^n(2M(\alpha ,A)+1)^n C_1(\alpha , A) a.\endaligned$$
Similarly we get
$$C_X ( \{ \psi -\varphi <-a\} )\leq 2^n(2M(\alpha ,A)+1)^n C_1(\alpha , A) a.$$
Combination of these inequalities yields
$$C_X(\{ |\varphi -\psi |>a \})\leq C(\alpha ,A)a.$$
Now we prove the promised generalization of Kolodziej stability theorem (Theorem 1.1 in [Ko5]).
\medskip
\noindent
{\bf 3.3. Theorem.} {\sl Let $\varphi,\psi\in\Cal E^-(X,\o)$ be such that $\sup\limits_X\varphi =\sup\limits_X\psi =0$ and $\o_\varphi^n, \o_\psi^n\in\Cal H(\alpha ,A)$. Then there exists $C(\alpha ,A)>0$ such that}
$$\sup\limits_X |\varphi -\psi |\leq C(\alpha ,A)[\int\limits_X ||\o_{\varphi}^n-\o_{\psi}^n||]^{\frac {\min (1,\frac \alpha n)} {2n+3+\frac {n+1}{1+\alpha}}}.$$ 
{\sl Proof.} Set 
$$a=[\int\limits_X ||\o_{\varphi}^n-\o_{\psi}^n||]^{ \frac 1 {2n+3+\frac {n+1}{1+\alpha}}}.$$
By Proposition 3.2 there exists $C_1(\alpha ,A)>0$ and $t\in\bold R$ such that $|t|\leq M(\alpha ,A)+1$ and
$$C_X(\{ |\varphi -\psi -t|>a \})\leq C_1(\alpha ,A)a.$$
Moreover, by Proposition 2.6 in [EGZ] there exists $C_2(\alpha ,A)>0$ such that
$$\aligned \sup\limits_X|\varphi-\psi -t|&\leq 2a + C_2(\alpha ,A)[C_X (\{|\varphi -\psi -t|>a \}) ]^{\frac \alpha n}\\
&\leq 2a+ C_2(\alpha ,A)[C_1(\alpha ,A)a]^{\frac \alpha n}\\
&\leq C_3(\alpha, A) a^{\min (1,\frac \alpha n)}.\endaligned$$
Moreover, since $\sup\limits_X\varphi =\sup\limits_X\psi =0$ we obtain $|t|\leq C_3(\alpha, A) a^{\min (1,\frac \alpha n)}$. Combination of these inequalities yields
$$\sup\limits_X|\varphi-\psi |\leq \sup\limits_X|\varphi-\psi -t|+|t|\leq 2C_3(\alpha, A) a^{\min (1,\frac \alpha n)}=C(\alpha ,A)[\int\limits_X ||\o_{\varphi}^n-\o_{\psi}^n||]^{ \frac {\min (1,\frac \alpha n)} {2n+3+\frac {n+1}{1+\alpha}}}.$$
{\bf 3.4. Corollary.} {\sl Let $\mu$ be a non-negative Radon measure on $X$ such that $\mu (B(z,r))\leq A r^{ 2n-2+\alpha }$ for all $B(z,r)\subset X$ ($A,\ \alpha >0$ are constants). Given $p>1, M>0, \epsilon >0$ and $f,g\in L^p (d\mu )$ with $||f||_{L^p(d\mu )}, ||g||_{L^p(d\mu )}\leq M$ and $\int\limits_X f d\mu = \int\limits_X g d\mu =1$. Assume that $\varphi,\psi\in\Cal E^-(X,\o)$ satisfy $\o_\varphi^n=fd\mu$, $\o_\psi^n=gd\mu$ and $\sup\limits_X\varphi =\sup\limits_X\psi =0$. Then there exists $C(\alpha ,A,M,\epsilon)>0$ such that}
$$\sup\limits_X |\varphi -\psi |\leq C(\alpha ,A, M, \epsilon )[\int\limits_X |f-g| d\mu ]^{ \frac {1} { 2n+3+\epsilon } }.$$ 
{\sl Proof.} By H\"older inequality we have 
$$\int\limits_K fd\mu \leq ||f||_{L^p(d\mu )}[\mu (K)]^{1-\frac 1 p}\leq M [\mu (K)]^{1-\frac 1 p},$$
$$\int\limits_K gd\mu \leq ||g||_{L^p(d\mu )}[\mu (K)]^{1-\frac 1 p}\leq M [\mu (K)]^{1-\frac 1 p},$$
for any Borel subset $K$ of $X$. By Proposition 2.7 we get $f d\mu, g d\mu \in \Cal H(\infty )$. Using Theorem 3.3 we can find $C(\alpha ,A,M,\epsilon)>0$ such that
$$\sup\limits_X |\varphi -\psi |\leq C(\alpha ,A, M, \epsilon )[\int\limits_X |f-g| d\mu ]^{ \frac {1} { 2n+3+\epsilon } }.$$ 
{\bf 4. Local estimates in Potential theory}
\medskip
\noindent
Let $\Omega$ be a bounded domain in $\bold R^n$ ($n\geq 2$). By SH$(\Omega )$ (resp SH$^-(\Omega )$) we denote the set of subharmonic (resp. negative subharmonic) functions on $\Omega$. For each $u\in SH (\Omega )$ and $\delta >0$ we denote
$$\tilde u_\delta (x)= \frac 1 {c_n\delta ^n} \int\limits_{B_\delta} u(x + y) dV_n (y),$$ 
$$u_\delta (x) =\sup\limits_{y\in B_\delta } u(x+y),$$
for $x\in\Omega_\delta =\{x\in\Omega: d(x,\partial\Omega )>\delta\}$. Here $B_\delta =\{x\in\bold R^n:\ |x|=(x_1^2+...+x_n^2)^{\frac 1 2}<\delta \}$ and $c_n$ is the volume of the unit ball $B_1$. We state some results which will be used in our main theorems. 
\medskip
\noindent
{\bf 4.1. Theorem.} {\sl Let $\mu$ be a non-negative Radon measure on $\Omega$ such that $\mu (B(z,r))\leq A r^{ n-2+\alpha }$ for all $B(z,r)\subset D\subset\subset\Omega$ ($A,\alpha >0$ are constants). Then for $K\subset\subset D$ and $\epsilon >0$ there exists $C(\alpha , A, K, \epsilon)$ such that
$$\int\limits_K [\tilde u_\delta- u] d\mu \leq C(\alpha ,A, K,\epsilon) \int\limits_{\bar D} \Delta u \ \delta^{\frac {\alpha -\epsilon}{1+\alpha }} ,$$
for all $u\in\text{SH}(\Omega)$, where $\Delta $ is the Laplace operator.}
\medskip
\noindent
{\sl Proof.} Since the change of radii of the balls does not affect the statement we can assume that $\Omega=B_4$, $D=B_3$, $K=B_1$ and $u$ is smooth on $B_4$. By [H\"o] we have  
$$u(x)=\int\limits_{B_2} G(x,z)\Delta u (z) +h (x),$$
where $G(x,y)$ is the fundamental solution of Laplace equation and $h$ is harmonic in $B_2$. By Fubini theorem we have
$$\aligned\int\limits_{B_1} [\tilde u_\delta (x)- u(x)] d\mu (x)=&\int\limits_{B_1} \frac 1 {c_n\delta ^n} \int\limits_{B_\delta} [u(x + y)-u(x)] dV_n (y) d\mu (x)\\
&\frac 1 {c_n\delta ^n}\int\limits_{B_1} \int\limits_{B_\delta} \int\limits_{B_2} [G(x+y,z)-G(x,z)]\Delta u (z) dV_n (y) d\mu (x)\\
&=\int\limits_{B_2}\Delta u(z) \frac 1 {c_n\delta ^n}\int\limits_{B_\delta} dV_n(y) \int\limits_{B_1} [G(x + y,z)-G(x,z)] d\mu (x) 
\endaligned.$$
Set 
$$F(y,z)= \int\limits_{B_1} [G(x + y,z)-G(x,z)] d\mu (x).$$
It is enough to prove that $F(y,z) \leq C(\alpha ,A,s) \delta^{\frac {\alpha -\epsilon} {1+\alpha } }$ for all $y\in B_\delta, z\in B_2$. We consider two cases:
\medskip
\noindent
Case 1: $n=2$. For $y\in B_\delta, z\in B_2$, $\delta <\frac 1 2$, we have
$$\aligned F(y,z)&=\int\limits_{B_1} [\ln |x + y-z|-\ln |x-z|] d\mu (x)\\
&=\int\limits_{ B_1\cap \{|x-z|\geq |y|^{\frac 1 {1+\alpha }} \} } \ln |1+\frac y {x-z}| d\mu (x)+\int\limits_{ B_1\cap \{|x-z|<|y|^{\frac 1 {1+\alpha }} \} }\ln |1+\frac y {x-z}| d\mu (x)\\
&\leq\int\limits_{ B_1\cap \{|x-z|\geq |y|^{\frac 1 {1+\alpha }} \} } \ln (1+|y|^{\frac {\alpha } {1+\alpha }}) d\mu (x)+\ln 4 \int\limits_{ B_1\cap \{|x-z|<|y|^{\frac 1 {1+\alpha }} \} } d\mu\\
&\ +\int\limits_{ B_1\cap \{|x-z|<|y|^{\frac 1 {1+\alpha }} \} }\ln \frac 1 {|x-z|} d\mu (x)\\
&\leq |y|^{\frac {\alpha} {1+\alpha }}\mu (B_1)+A |y|^{\frac {\alpha} {1+\alpha }}\ln 4+|y|^{\frac {\alpha -\epsilon} {1+\alpha} } \int\limits_{\{|x-z|<|y|^{\frac 1 {1+\alpha }} \} }\frac {1} { |x-z|^{ \alpha -\epsilon } }\ln \frac 1 {|x-z|} d\mu (x)\\
&\leq A(1+\ln 4 )|y|^{\frac {\alpha} {1+\alpha }}+|y|^{\frac {\alpha -\epsilon} {1+\alpha} } C_1(\alpha ,\epsilon )\int\limits_{\{|x-z|<1 \} }\frac {d\mu (x)} { |x-z|^{ \alpha-\frac {\epsilon} 2 } } \\
&\leq A(1+\ln 4 )|y|^{\frac {\alpha} {1+\alpha }}+C_1(\alpha ,\epsilon )|y|^{\frac {\alpha -\epsilon} {1+\alpha} }\sum\limits_{j=0}^\infty\int\limits_{\{ 2^{-j-1}\leq |x-z|< 2^{-j} \}} \frac {d\mu (x)} { |x-z|^{ \alpha-\frac {\epsilon} 2 } } \\
&\leq A(1+\ln 4 )|y|^{\frac {\alpha} {1+\alpha }}+C_1(\alpha ,\epsilon )|y|^{\frac {\alpha -\epsilon} {1+\alpha} }A \sum\limits_{j=0}^\infty  2^{ (j+1)(\alpha -\frac \epsilon 2)-j\alpha }\\
&\leq C(\alpha, A, \epsilon ) |y|^{\frac {\alpha -\epsilon} {1+\alpha }}\leq C(\alpha ,A, \epsilon ) \delta^{\frac {\alpha -\epsilon} {1+\alpha }}.
\endaligned$$
Case 2: $n\geq 3$. Similarly for $y\in B_\delta, z\in B_2$, $\delta <\frac 1 2$, we have
$$\aligned F(y,z)&=\int\limits_{B_1} [-\frac 1 {|x + y-z|^{n-2}}+\frac 1{|x-z|^{n-2}}] d\mu (x)\\
&=\int\limits_{ B_1\cap \{|x-z|\geq |y|^{\frac 1 {1+\alpha }} \} } \frac {|x + y-z|^{n-2}-|x-z|^{n-2}} {|x + y-z|^{n-2}|x-z|^{n-2} } d\mu (x)+\int\limits_{\{|x-z|<|y|^{\frac 1 {1+\alpha }} \} }\frac {d\mu (x)} {|x-z|^{n-2}}\\
&\leq C_2(\alpha ) |y|^{\frac {\alpha } {1+\alpha }}\int\limits_{ B_1\cap \{|x-z|\geq |y|^{\frac 1 {1+\alpha }} \} } d\mu (x)+|y|^{\frac {\alpha -\epsilon} {1+\alpha } }\int\limits_{\{|x-z|<|y|^{\frac 1 {1+\alpha }} \} } \frac {d\mu (x)} {|x-z|^{n-2+\alpha-\epsilon }}\\
&\leq A C_2(\alpha ) |y|^{\frac {\alpha } {1+\alpha }}+|y|^{\frac {\alpha -\epsilon} {1+\alpha } }\int\limits_{\{ |x-z|<1 \} } \frac {d\mu (x)} {|x-z|^{n-2+\alpha-\epsilon }}\\
&\leq C(\alpha, A,\epsilon) |y|^{\frac {\alpha -\epsilon} {1+\alpha }}\leq C(\alpha ,A, \epsilon) \delta^{\frac {\alpha -\epsilon}{1+\alpha }},
\endaligned$$
{\bf 4.2. Theorem.} {\sl Let $\mu$ be a non-negative Radon measure on $\Omega$ such that $\mu (B(z,r))\leq A r^{ n-2+\alpha }$ for all $B(z,r)\subset D\subset\subset\Omega$ ($A,\alpha >0$ are constants). Then for $K\subset\subset D$ and $\epsilon>0$ there exists $C(\alpha , A, K, \epsilon)$ such that}
$$\int\limits_K [u_\delta- u] d\mu \leq C(\alpha ,A, K,\epsilon) ||u||_{L^\infty (\Omega )} \ \delta^{\frac {\alpha -\epsilon}{2(1+\alpha )}},$$
{\sl for all $u\in\text{SH}\cap L^\infty (\Omega)$.}
\medskip
\noindent
We need a well-known lemma:
\medskip
\noindent
{\bf 4.3. Lemma.} {\sl Let $u\in\text{SH}\cap L^\infty (\Omega)$. Then}
$$|\tilde u_\delta (x)-\tilde u_\delta (y)|\leq \frac {||u||_{L^\infty (\Omega)}|x-y|}{\delta },$$
for all $x,y\in\Omega_\delta$.
\medskip
\noindent
{\sl Proof of Theorem 4.2.} By Lemma 4.3 we have
$$u_\delta (x)=\sup\limits_{y\in B_\delta } u(x+y)\leq \sup\limits_{y\in B_\delta } \tilde u_{\delta^{\frac 1 2}}(x+y)\leq \tilde u_{\delta^{\frac 1 2}}(x)+\delta ^{\frac 1 2} ||u||_{L^\infty (\Omega )}.$$
By Theorem 4.1 we get
$$\aligned\int\limits_K [u_\delta- u] d\mu&\leq\int\limits_K [\tilde u_{\delta ^{\frac 1 2}}- u]d\mu +||u||_{L^\infty (\Omega)} \mu (K)\delta ^{\frac 1 2}\\
&\leq C(\alpha ,A, K,\epsilon) ||u||_{L^\infty (\Omega ) } \ \delta^{\frac {\alpha -\epsilon}{2(1+\alpha )}}
.\endaligned$$
Next we state a well-known result is a direct consequence of the Jensen formula (see [AG])
\medskip
\noindent
{\bf 4.4. Proposition.} {\sl Let $u\in\text{SH}(B_2)$ be such that $|u(x)-u(y)|\leq A|x-y|^{\alpha}$ for all $x,y\in B_2$. Then there exists $C(\alpha, A)>0$ such that
$$\int\limits_{B(x,r)}\Delta u\leq C(\alpha, A) r^{ n-2+\alpha },$$ 
for all $B(x,r)\subset B_1$.}
\vskip 0.5 cm
\noindent
{\bf 5. Main results}
\medskip
\noindent
{\sl Proof of Theorem A.} We use the same scheme as the proof of Theorem 2.1 in [Ko5]. From Corollary 3.4 and from Theorem 4.2 we can replace $\omega^n$ by $d \mu$. This implies that $u$ is H\"older continuous with the H\"older exponent dependent on $\alpha $, $A$, $p$, $X$ and $||f||_{L^p(d\mu )}$.
\medskip
\noindent
{\sl Proof of Corollary B.} It follows from Proposition 4.4 and Theorem A.
\medskip
\noindent
{\sl Proof of Corollary C.} Direct application of Theorem A.
\vskip 0,5cm
\noindent
\centerline{\bf REFERENCES}
\vskip 0,5cm
\noindent
[AG] D. H. Armitage and S. J. Gardiner, Classical potential theory, Springer Monogr. Math., Springer-Verlag, London, 2001.

\noindent
[BT1] E. Bedford and B. A. Taylor, The Dirichlet problem for the complex Monge-Amp\`{e}re operator, Invent. Math. 37 (1976), 1-44.

\noindent
[BT2] E. Bedford and B. A. Taylor, A new capacity for plurisubharmonic functions, Acta Math. 149 (1982), 1-40.

\noindent
[Ce1] U. Cegrell, Pluricomplex energy, Acta Math. 180 (1998), 187-217.

\noindent
[Ce2] U. Cegrell, The general definition of the complex Monge-Amp\`{e}re operator, Ann. Inst. Fourier (Grenoble) 54 (2004), 159-179.

\noindent
[CK] U. Cegrell, S. Ko{\l}odziej, The equation of complex Monge-Amp\`ere type and stability of solutions, Math. Ann. 334 (2006), 713-729.

\noindent
[CGZ] D. Coman, V. Guedj and A. Zeriahi, Domains of definition of Monge-Amp\`{e}re operators on compact K\"ahler manifolds, Math. Zeit. 259 (2008), 393-418.

\noindent
[De1] J. P. Demailly, Mesures de Monge-Amp\`{e}re et mesures pluriharmoniques, Math. Zeit. 194 (1987), 519-564.

\noindent
[De2] J. P. Demailly, Monger-Amp\`{e}re operators, Lelong numbers and intersection theory, Complex Analysis and Geometry, Univ. Series in Math., Plenum Press, New-York, 1993.

\noindent
[De3] J. P. Demailly, Complex analytic and differential geometry, self published e-book, 1997. 

\noindent
[Di1] S. Dinew, Cegrell classes on compact K\"ahler manifolds, Ann. Polon. Math. 91 (2007), 179-195.

\noindent
[Di2] S. Dinew, An inequality for mixed Monge-Amp\`{e}re measures, Math. Zeit. 262 (2009), 1-15.

\noindent
[Di3] S. Dinew, Uniqueness in $\Cal E(X,\o )$, J. Funct. Anal. 256 (2009), 2113-2122.

\noindent
[DH] S. Dinew and P. H. Hiep, Convergence in capacity on compact K\"ahler manifolds, Preprint (2009) (http://arxiv.org). 

\noindent
[DZ] S. Dinew and Z. Zhang, Stability of Bounded Solutions for Degenerate Complex Monge-Amp\`{e}re equations, Preprint (2008) (http://arxiv.org)

\noindent
[DNS] T. C. Dinh, V. A. Nguyen and N. Sibony, Exponential estimates for plurisubharmonic functions and
stochastic dynamics, Preprint (2008) (http://arxiv.org)

\noindent
[DS] T. C. Dinh, N. Sibony, Distribution des valeurs de transformations mեromorphes et applications, Comment. Math. Helv., 81 (2006), no. 1, 221-258.

\noindent
[EGZ] P. Eyssidieux, V. Guedj and A. Zeriahi, Singular K\"ahler Einstein metrics, J. Amer. Math. Soc. 22 (2009), 607-639. 

\noindent
[GZ1] V. Guedj and A. Zeriahi, Intrinsic capacities on compact K\"ahler manifolds, J. Geom. Anal. 15 (2005), no. 4, 607-639.

\noindent
[GZ2] V. Guedj and A. Zeriahi, The weighted Monge-Amp\`{e}re energy of quasiplurisubharmonic functions, J. Funct. Anal. 250 (2007), 442-482.

\noindent
[H] P. H. Hiep, On the convergence in capacity on compact K\"ahler manifolds and its applications, Proc. Amer. Math. Soc. 136 (2008), 2007-2018.

\noindent
[H\"o] L. H\"ormander, Notions of Convexity, Progess in Mathematics 127, Birkh\"auser, Boston (1994).

\noindent
\noindent
[Ko1] S. Ko{\l}odziej, The Monge-Amp\`{e}re equation, Acta Math. 180 (1998), 69-117.

\noindent
[Ko2] S. Ko{\l}odziej, The Monge-Amp\`{e}re equation on compact K\"ahler manifolds, Indiana Univ. Math. J. 52 (2003), 667-686.

\noindent
[Ko3] S.Ko{\l}odziej, The set of measures given by bounded solutions of the complex Monge-Amp\`ere equation on compact K\"ahler manifolds, J. London Math. Soc. 72 (2) (2005), 225-238.

\noindent
[Ko4] S. Ko{\l}odziej, The complex Monge-Amp\`{e}re equation and pluripotential theory, Mem. Amer. Math. Soc. 178 (2005).

\noindent
[Ko5] S. Ko{\l}odziej, H\"older continuity of solutions to the complex Monge-Amp\`{e}re equation with the right-hand side in $L^p$: the case of compact K\"ahler manifolds, Math. Ann. 342 (2008), 379-386.

\noindent
[KoTi] S. Ko{\l}odziej and G. Tian, A uniform $L^\infty$ estimate for complex Monge-Amp\`{e}re equations, Math. Ann. 342 (2008), 773-787.

\noindent
[Si1] J. Siciak, On some extremal functions and their applications in the theory of analytic functions of several complex variables, Trans. Amer. Math. Soc. 105 (1962) 322ֳ57.

\noindent
[Si2] J. Siciak, Extremal plurisubharmonic functions and capacities in $\bold C^n$, Sophia Univ., Tokyo (1982).

\noindent
[Ze1] A.Zeriahi, The size of plurisubharmonic lemniscates in terms of Hausdorff-Riesz measures and capacities, Proc. London Math. Soc. 89 (2004), no. 1, 104-122.

\noindent
[Ze2] A. Zeriahi, A minimum Principle for Plurisubharmonic functions, Indiana Univ. Math. J. 56 (2007), 2671-2696.

\noindent
[Yau] S. T. Yau, On the Ricci curvature of a compact K\"ahler manifold and the complex Monge-Amp\`{e}re equation, Commun. Pure Appl. Math. 31 (1978), 339-411.
\medskip
\noindent
Department of Mathematics
\medskip
\noindent
University of Education (Dai hoc Su Pham  Ha Noi)
\medskip
\noindent
CauGiay, Hanoi, Vietnam
\medskip 
\noindent
E-mail: phhiep$_-$vn$\@$yahoo.com
\end